\documentclass[reqno,12pt,a4paper]{amsart}
\usepackage{pgfplots}
\usepackage[english]{babel}
 \pdfoutput=1
\usepackage{amsmath,amsthm,amssymb,amsfonts, mathdots}
\usepackage{bbm}
\usepackage{scrtime, cancel}
\usepackage{enumitem, color, comment}
  
\renewcommand{\le}{\leqslant}
\renewcommand{\leq}{\leqslant}
\renewcommand{\Re}{\mathrm{Re}}

\renewcommand{\ge}{\geqslant}

\definecolor{pinegreen}{rgb}{0.0, 0.47, 0.44}
\definecolor{brilliantrose}{rgb}{1.0, 0.33, 0.64}
\usepackage{mathtools}
\usepackage{hyperref}

\usepackage{pdfsync,verbatim}
\usepackage{color}
\mathtoolsset{showonlyrefs}

\numberwithin{equation}{section}   

\allowdisplaybreaks
\newtheorem{theorem}{Theorem}[section]
\newtheorem{lemma}[theorem]{Lemma}
\newtheorem{proposition}[theorem]{Proposition}
\newtheorem{corollary}[theorem]{Corollary}

\theoremstyle{definition}
\newtheorem{remark}[theorem]{Remark}

\definecolor{pinegreen}{rgb}{0.0, 0.47, 0.44}
\newcommand{\R}{\mathbb{R}}

\newcount\dotcnt\newdimen\deltay
\def\Ddot#1#2(#3,#4,#5,#6){\deltay=#6\setbox1=\hbox to0pt{\smash{\dotcnt=1
\kern#3\loop\raise\dotcnt\deltay\hbox to0pt{\hss#2}\kern#5\ifnum\dotcnt<#1
\advance\dotcnt 1\repeat}\hss}\setbox2=\vtop{\box1}\ht2=#4\box2}

\makeatletter
\def\author@andify{%
  \nxandlist {\unskip ,\penalty-1 \space\ignorespaces}%
    {\unskip {} \@@and~}%
    {\unskip \penalty-2 \space \@@and~}%
}
\makeatother

 \title[Spectrum  of the inverse Ornstein--Uhlenbeck 
 operator ]{Spectral and mapping 
 properties \\of the inverse Ornstein--Uhlenbeck operator 
 }

 \author[T.\ Bruno]{Tommaso Bruno}
\address{Dipartimento di Matematica, Universit\`a degli Studi di Genova\\Via Dodecaneso 35 \\16146 Genova\\ Italy}
\email{tommaso.bruno@unige.it}
\author[V.\ Casarino]{Valentina Casarino}
\address{DTG, Universit\`a degli Studi di Padova\\ Stradella san Nicola 3 \\I-36100 Vicenza \\ Italy}
\email{valentina.casarino@unipd.it}
\author[P.\ Ciatti]{Paolo Ciatti}
\address{Dipartimento di Matematica ``Tullio Levi Civita'', Universit\`a degli Studi di Padova\\Via Trieste, 63, 35131 Padova,  \\ Italy}
\email{paolo.ciatti@unipd.it}
\author[P.\ Sj\"ogren]{Peter Sj\"ogren}
\address{Mathematical Sciences,  University of Gothenburg and  Mathematical Sciences,
Chalmers University of Technology  \\ SE - 412 96 G\"oteborg, Sweden}
\email{peters@chalmers.se}

\thanks{The first 
three authors are members 
of the Gruppo Nazionale per l'Analisi Matematica, la Probabilità e le loro Applicazioni (GNAMPA)
of the Istituto Nazionale di Alta Matematica (INdAM).
Part of this research was carried out while the fourth author was GNAMPA Professore Visitatore
at the University of Padova in Italy. The fourth author also profited from a grant from the Royal Society
of Arts and Sciences in Gothenburg and another from the Royal Swedish Academy of Sciences. He
is also grateful to the University of Padova for its hospitality during several visits.}


\date
\today

\begin{document}

\begin{abstract}
We investigate the spectral and boundedness properties 
of the inverse Ornstein--Uhlenbeck operator $\mathcal{A}^{Q,B}$
and its associated semigroup  in $L^p$ spaces with respect to the inverse Gaussian measure. We establish  strong continuity and contractivity of this semigroup in $L^p$ for $p \in [1, \infty)$, and obtain a complete characterization 
of the  spectrum of $\mathcal{A}^{Q,B}$.
For $p \in (1, \infty)$, the spectrum is independent of $p$, 
consists   entirely of eigenvalues 
and exhibits a spectral gap.     
For $p=1$, the spectrum is instead  the left half-plane. 
Furthermore, we  discuss ultracontractivity  
and  hypercontractivity properties 
of  the semigroup.                 
\end{abstract}

\maketitle

\section{Introduction}\label{Introduction}
Throughout this paper, we consider two real $n\times n$ matrices $Q$ and $B$ fulfilling the following two conditions:
\begin{enumerate}
 \item[(H1)] $Q$ is symmetric and positive definite;
 \item[(H2)] all the eigenvalues of $B$ have negative real parts.
\end{enumerate}
For each $t\in (0,+\infty]$, we introduce the covariance matrices
\begin{equation}\label{defQt}
Q_t=\int_0^t e^{sB}Qe^{sB^*}ds \,.
\end{equation}
Each $Q_t$ is well defined, symmetric and positive definite. 
Thus $Q_t$ increases in the sense of quadratic forms and converges to  
$Q_\infty$.   
\\                               
We consider $\R^n$ endowed with the inverse Gaussian measure
\begin{align}\label{def: gaussian_measure}
d{\gamma_{-\infty}}(x)=
\exp\big({\tfrac12 \langle Q_\infty^{-1} x, x\rangle}\big)\,dx \,.
\end{align}
Our aim is to study spectral properties of the inverse Ornstein--Uhlenbeck operator, given by
\begin{equation}\label{def:opA}
\mathcal A^{Q,B}=
\tfrac12\,
\mathrm{tr}
\big( Q\nabla^2 \big)-\langle Bx, \nabla \rangle,
\end{equation}
in $L^p(\gamma_{-\infty})$. Here the symbols $\nabla$ and $\nabla^2$ denote the gradient and the Hessian, respectively. 

The semigroup generated by $\mathcal A^{Q,B}$ will be denoted by the symbol    
$$
\mathcal H_t^{UO} = e^{t \mathcal A^{Q,B}}\,, \qquad t>0.
$$
As in \cite{BCCS}, 
we use the superscript ``$UO$'' to highlight the relationship with the Ornstein--Uhlenbeck setting, 
where the superscript ``$OU$'' is used. 
We recall that the classical Ornstein--Uhlenbeck semigroup $ \mathcal H_t^{OU}$ 
is  generated by the Ornstein--Uhlenbeck operator 
\begin{equation}\label{generalizedOU}
\mathcal L^{Q,B}=
\tfrac12\,
\mathrm{tr}
\big( Q\nabla^2 \big)+\langle Bx, \nabla \rangle,
\end{equation} 
where $Q$ and $B$ are the same as in \eqref{def:opA}. 
For  extensive information on recent harmonic analysis  
concerning the Ornstein--Uhlenbeck semigroup, 
we refer the reader to \cite{CCS2, CCS3, CCS5,  CCS8, CCS9} 
and the references therein.

The operator $\mathcal A^{Q,B}$ has been introduced by Salogni in \cite{Salogni}
 in the standard case $Q=-B=I_n$, where $I_n$ denotes the identity matrix;
 there, the author studied $L^p$ boundedness properties of the maximal operator for the associated semigroup.
Under the same assumption $Q=-B=I_n$, 
the first and the fourth authors of this paper studied boundedness properties 
of some Riesz transforms associated to $\mathcal H_t^{UO}$ (see \cite{Bruno2, BS}).
Other results related to $\mathcal A^{I,-I}$ may be found in \cite{AB, ABFR, BMRM, BCLC, Mazzitelli}.
More recently, we introduced a more general version of $\mathcal A^{Q,B}$, with $Q$ and $B$ only fulfilling (H1) and (H2), 
and proved the weak type (1,1) for the associated maximal function \cite{BCCS}. 
Notice that the measure used in \cite{BCCS} differs from $d\gamma_{-\infty}$ by an irrelevant constant factor.

The study of the inverse Ornstein--Uhlenbeck operator presents several challenges 
that do not arise in the classical Ornstein--Uhlenbeck setting. 
Most notably, although $\gamma_{-\infty}$ is an invariant measure for $(\mathcal{H}_t^{UO})_{t>0}$, such semigroup does not admit an invariant \emph{probability} measure; and
unlike the symmetric inverse setting $Q=-B= I_{n}$ considered in \cite{Salogni,Bruno2,BS},
neither the operator nor the associated semigroup is symmetric in general.
Furthermore, the  reference measure $\gamma_{-\infty}$ is  locally doubling on admissible balls, see~\cite{Bruno2},
 but lacks the global doubling property, 
and it exhibits superexponential growth at infinity. 
This contrast between the global growth of the measure 
and its local Euclidean behavior 
strongly affects the boundedness and contractivity of the semigroup.

The main purpose of this paper is to provide a comprehensive analysis 
of the spectral and boundedness properties 
of $\mathcal{A}^{Q,B}$ and its associated semigroup. 
First, we establish the strong continuity and the contractivity of $(\mathcal{H}_t^{UO})_{t>0}$
in $L^p(\gamma_{-\infty})$ for all $p \in [1, \infty)$. 
Building on these fundamental properties, 
we obtain a complete characterization of the spectrum of the generator $\mathcal{A}^{Q,B}$ in $L^p(\gamma_{-\infty})$. 
We prove that for $p \in (1, \infty)$ 
the spectrum is independent of $p$,  consists 
entirely of eigenvalues determined by those of the matrix $B+QQ_\infty^{-1}$, and
 exhibits a spectral gap. In  contrast, for $p=1$ the spectrum is a left half-plane.

Furthermore, we investigate the smoothing behavior of the semigroup.
We show that $(\mathcal H_t^{UO})_{t>0}$ is ultracontractive
with respect to both the Lebesgue measure and the inverse Gaussian measure.
This is in contrast with the Ornstein--Uhlenbeck semigroup. Finally, we prove hypercontractivity with respect to
$\gamma_{-\infty}$.

\medskip

\subsection{Scheme of the paper} 
Section \ref{section:measure} collects some preliminary facts about the inverse Gaussian measure, the inverse Ornstein--Uhlenbeck semigroup and the associated inverse  Mehler kernel. 
In Section \ref{section:cont}, 
we prove the strong continuity and the $L^p$-contractivity of the semigroup. 
Section \ref{section:spectrum} is devoted 
to the description of the spectrum of the generator $\mathcal{A}^{Q,B}$.

 We also discuss the orthogonality or the failure of orthogonality of some generalized eigenspaces 
in the inverse setting. 
In Section \ref{section:ultra}, 
we establish the ultracontractivity of the semigroup and discuss a first counterexample concerning its short-time behavior. 
Finally, Section \ref{sec:hyper}  provides a detailed analysis of the hypercontractivity of $(\mathcal{H}_t^{UO})_{t>0}$.

\subsection{Notation}
We shall adopt the ``variable constant convention", according to which positive constants are denoted by symbols such as $c>0$ or $C<\infty$, and these are not necessarily equal at different occurrences.  All constants will  depend only on $n$,  $Q$ and  $B$, unless otherwise explicitly stated.
                    For any two non-negative quantities $a$ and $b$ we write  $a\lesssim b$ (or $b\gtrsim a$) if $a\leq Cb$ for some $C$;
                     moreover $a\simeq 
b$ stands for the conjunction of $a\lesssim b$ and $b\lesssim a$.
By $\mathbb N$ we mean $\{0,1,\dots\}$. 

\section{The inverse Gaussian measure}\label{section:measure}
As  is well known,  the Gaussian probability measure
\begin{equation}\label{def:meas_gauss}
d\gamma_\infty (x)=
(2\pi)^{-\frac{n}{2}}
(\det\, Q_\infty)^{-\frac{1}{2} }
\exp\big({-\tfrac12 \langle Q_\infty^{-1} x, x\rangle}\big)\,dx
\end{equation}
is an invariant (probability) measure 
for the semigroup generated by $\mathcal L^{Q,B}$, that is,
\begin{align*}
 \int_{\R^{n}} \mathcal H_t^{OU} f(x)\, d\gamma_\infty(x)=\int_{\R^{n}} f(x)\, d\gamma_\infty(x)\qquad \text{for all\; $t>0$ \;and \;$f\in L^1(\gamma_\infty )$.}
\end{align*}
In fact, in \cite[Section 11.2.3  page 313]{DaPratoZ}
it is proved that the existence of an invariant probability
measure for the semigroup generated by an operator $\mathcal L^{Q,B}$ 
(for an arbitrary $B$)
defined as in  \eqref{generalizedOU}  is equivalent to 
the spectrum of the matrix  $B$  
 being contained in the open left half-plane (we refer also to~\cite{Metafune} or~\cite{Lunardi}).
Since  the spectrum of   $-B$   is contained in the open right half-plane, $\mathcal H_t^{UO}$ has no
invariant {probability} measure.

We shall see in Proposition~\ref{pr:contr_L1_pd}, however, that $d\gamma_{-\infty}$ can be seen as an invariant measure for $\mathcal H_t^{UO}$, despite not being finite. An additional reason for which $d\gamma_{-\infty}$ is a natural measure in the inverse Ornstein--Uhlenbeck setting is given by the following lemma.

\begin{lemma}
The semigroup  $e^{t\mathcal{A}^{Q,B}}$  is symmetric on $L^{2}(\gamma_{-\infty})$ if and only if $QB^* = BQ$.
\end{lemma}
\begin{proof}
It has been observed  in \cite[Remark 2.3]{BCCS} that the semigroup $e^{t\mathcal{A}^{Q,B}}$ is symmetric on $L^{2}(\gamma_{-\infty})$ if and only if $e^{t\mathcal{L}^{Q,B}}$ is symmetric on $L^{2}(\gamma_{\infty})$. Since the latter property 
holds if and only if $QB^* = BQ$ by~\cite[Theorem 2.4]{CMG} (see also \cite[Lemma 2.1]{MPRS}),
the lemma follows.
\end{proof}

\medskip

 \medskip

We recall from \cite{BCCS} the following facts. To simplify the notation, let us first set
\begin{equation}\label{def:cn}
c_n= (2\pi)^{-\frac{n}{2}} (\det \, Q_\infty)^{-\frac{1}{2} }.
\end{equation}
Denoting by $R(x)$ the quadratic form
\begin{equation}\label{quadratic-form}
R(x) ={\frac12 \left\langle Q_\infty^{-1}x ,x  \right\rangle}, \qquad\text{$x\in\R^n$},
\end{equation}
we have
$
d\gamma_\infty (x)=c_n \, e^{-R(x) }\,dx$,
 as well as 
 $d{\gamma_{-\infty}}(x) = e^{R(x) }\,dx$,
and thus
\begin{equation}\label{eq:rel_meas}
d\gamma_\infty (x)=c_n e^{-2R(x)}\, d{\gamma_{-\infty}}(x). 
\end{equation}

 \medskip

\begin{enumerate}
 \item 
 For each $f\in L^1(\gamma_{-\infty})$ and all $t>0$ one has
\begin{align}
e^{t \mathcal A^{Q,B}} f(x) &=
 \int
K_t^{UO}
(x,u)\,
f(u)\,
 d\gamma_{-\infty}(u)
  \,, \qquad x\in\R^n,
\end{align}
 where for $x,u\in\R^n$ and $t>0$ the inverse Mehler kernel $K_t^{UO}$ (with respect to $d\gamma_{-\infty}$)   is given by
\begin{align*}   
  K_t^{UO}(x,u) &=
(2\pi)^{-\frac n2} ( \det Q_t)^{-\frac12}\,  \\
& \qquad  \times e^{t\,\mathrm{tr}\, B}\,  e^{-R(x)}
   \, \exp\left[ -\frac12
  \langle (Q_t^{-1} -Q_\infty^{-1})(x-D_t\, u), \,x-D_t\, u \rangle \right]
  \end{align*}
 (see \cite[Proposition 2.2]{BCCS} for a proof).
 Here $D_t $ denotes a  one-parameter group of matrices, given by
\begin{equation}\label{def:Dtx}
D_t := Q_\infty \,
 e^{-tB^*}\, Q_\infty^{-1},  \quad \; t \in \R.
\end{equation}
 We also recall \cite[formula (4.3)]{CCS2}, that is,
\begin{align}
&\frac{\partial}{\partial t} \,
R\big( D_t x \big)
\simeq
\big| 
D_t x\big|^2,
 \label{vel-4}
\end{align}
for all  $x$
in $\R^n$ and $t\in\R$. 
This implies that for $|t|<1$
\begin{equation} \label{close}
|D_t \,x -x| \simeq |t| |x|.
\end{equation}     
 \item 
 Turning  back to the Ornstein--Uhlenbeck setting, one has for each $f\in L^1(\gamma_\infty)$ and all $t>0$ 
\begin{align}                     
e^{t \mathcal L^{Q,B}} f(x) &=
 \int
K_t^{OU}
(x,u)\,
f(u)\,
 d\gamma_\infty(u)
  \,, \qquad x\in\R^n,
\end{align}
 where for $t>0$ the Mehler kernel $K_t^{OU}$ (with respect to $d\gamma_\infty$)   is given by
\begin{align}\label{mehler}   
K_t^{OU} (x,u)
&=
\Big(
\frac{\det \, Q_\infty}{\det \, Q_t}
\Big)^{\frac{1}{2} }
e^{ R(x)}
\exp \bigg[
{-\frac12 
\left\langle (
Q_t^{-1}-Q_\infty^{-1}) (u-D_t x) \,,\, u-D_t x\right\rangle}\bigg]
\,
\end{align}
for $x,u\in\R^n$ (see~\cite[(2.6)]{CCS2}).
\item 
The relationship between the integral kernels  $K_t^{OU}$ and $K_t^{UO} $
is given by 
\begin{equation}\label{duenuclei}
 K_t^{UO}(x,u) =
c_n\,  e^{-R(x)} \, e^{-R(u)}\, e^{t\,\mathrm{tr}\, B}\,K_t^{OU}(u,x),
\end{equation}
holding for all $t>0$ and $x,u\in\mathbb R^n$ (see \cite[(2.5)]{BCCS}).

\item
In 
 \cite[Proposition 2.5]{BCCS} we proved the following
 Kolmogorov-type formula in the inverse Gaussian setting. 
 For all $f\in \mathcal C_b(\R^n)$ and $x\in\R^n$
 one has
 \begin{equation}\label{kolmogo}
\mathcal H_t^{UO}f(x)=(2\pi)^{-\frac{n}{2}}({\det} \, Q_t)^{-\frac{1}{2} }\int_{\mathbb{R}^{n}}f(e^{-tB}(u+x))  \, e^{-\frac12 \langle Q_t^{-1}u,u\rangle}\, du.
\end{equation}
\item
The semigroup  $(\mathcal H_t^{UO})_{t>0}$ is conservative,  that is, 
\begin{equation}\label{cons_eq}
 \mathcal H_t^{UO} \, \mathbf{1}=\mathbf{1}  \qquad \text{ for all $t>0$,}
\end{equation}
where $\mathbf{1}$ denotes
the constant 
function.
\item
 If $0 < t\le 1$, one has
\begin{equation}\label{litet}
 \frac{ e^{-R( x)}}{t^{n/2}}  \exp\bigg(-C\,\frac{|D_{-t}\,x- u |^2}t\bigg)
 \lesssim   K_t^{UO}(x,u)
 \lesssim  \frac{ e^{-R( x)}}{t^{n/2}} \exp\bigg(-c\,\frac{|D_{-t}\,x- u |^2}t\bigg)
\end{equation}
for all $(x, u) \in\mathbb R^n \times \mathbb R^n$ (see \cite[Proposition 2.4]{BCCS}).
\end{enumerate}

\section{On the  continuity and  contractivity of $(\mathcal H_t^{UO})_{t>0}$}\label{section:cont}

We start by considering some boundedness properties of $\big( \mathcal H_t^{UO}\big)_{t>0}$.

\begin{proposition}\label{pr:contr_L1_pd}
The measure $\gamma_{-\infty}$ is an invariant (infinite) measure for $\mathcal{H}_{t}^{UO}$, in the sense that 
\[
\int_{\R^{n}} \mathcal{H}_{t}^{UO} f\, d\gamma_{-\infty} =\int_{\R^{n}}f\, d\gamma_{-\infty}
\]
for all  $f \in L^1(\gamma_{-\infty})$ and $t>0$. Moreover, for each $t>0$
 and for all $1\le p\le \infty$
the operator $\mathcal H_t^{UO}$ is a contraction on $L^p(\gamma_{-\infty})$.
\end{proposition}

\begin{proof}
Let $f\in L^1(\gamma_{-\infty})$.
We may assume, without loss of generality, that $f$ is nonnegative.
Then
\begin{align*}
\int_{\R^{n}} \mathcal H_t^{UO} f(x) \, d\gamma_{-\infty}(x)&=
\int_{\R^{n}}\int_{\R^{n}} K_t^{UO} (x,u) \,f(u)\, d\gamma_{-\infty}(u)\, d\gamma_{-\infty}(x)\\
&=
\int_{\R^{n}}\int_{\R^{n}} K_t^{UO} (x,u) \,d\gamma_{-\infty}(x)\, f(u)\, d\gamma_{-\infty}(u).
\end{align*}
Setting 
\begin{align*}
 I_t (u)&:=
\int_{\R^{n}} K_t^{UO} (x,u) \,d\gamma_{-\infty}(x),\end{align*}
it will be enough to verify that
\begin{equation}\label{th_Itu}
I_t(u) = 1 
\end{equation}
for all $u\in\R^n$.

In the light of \eqref{duenuclei},  \eqref{def: gaussian_measure} and \eqref{mehler},
we observe that
\begin{align}
 I_t (u)
 &=
c_{n}\,\int_{\R^{n}}    e^{-R(x)} \, e^{-R(u)}\, e^{t\,\mathrm{tr}\, B}\,K_t^{OU}(u,x) \,d\gamma_{-\infty}(x)
\notag\\
 &=
c_{n}\,\int_{\R^{n}}   e^{-R(u)}\, e^{t\,\mathrm{tr} \,B}\,K_t^{OU}(u,x) \,dx
\label{f:intermediate}\\
 &=
 (2\pi)^{-\frac n2}\,e^{t\,\mathrm{tr} \,B}\, \big(
{\det \, Q_t}
\big)^{-\frac12 }
\int_{\R^{n}} \exp \bigg[
{-\frac12
\left\langle \left(
Q_t^{-1}-Q_\infty^{-1}\right) (x-D_t \,u) \,, x-D_t\, u\right\rangle}\bigg] \,
dx
\notag\\
 &=
 (2\pi)^{-\frac n2}\,e^{t\,\mathrm{tr} \,B}\, \big(
{\det \, Q_t}
\big)^{-\frac12 }
\int_{\R^{n}} \exp \bigg[
{-\frac12
\left\langle \left(
Q_t^{-1}-Q_\infty^{-1}\right) x \,, x\right\rangle}\bigg] \,
dx
\notag\\
 &=
 (2\pi)^{-\frac n2}\,e^{t\,\mathrm{tr} \,B}\, \big(
{\det \, Q_t}
\big)^{-\frac12 }
\big(\det 
\left(
Q_t^{-1}-Q_\infty^{-1}\right)\big)^{-\frac12 }
\int_{\R^{n}} \exp \left[
{-\frac12
\left\langle x' \,, x'\right\rangle}\right] \,
dx'.\label{conto_contrvv}
\end{align}
We know from
\cite[formulae (2.5) and (2.2)]{CCS2} that
 \begin{equation*}
Q_\infty- Q_t=
e^{tB} Q_\infty e^{tB^*}
\end{equation*}
(see also \cite[formula (2.1)]{MPRS}),
and  
\begin{equation*} 
Q_t ^{-1} -  Q_\infty^{-1} = Q_t^{-1}( Q_\infty -Q_t )Q_\infty^{-1}.\end{equation*}
It follows that
\begin{equation*}
Q_t ^{-1} -  Q_\infty^{-1} = Q_t^{-1}\, e^{tB} \, Q_\infty\, e^{tB^*} \,Q_\infty^{-1},
\end{equation*}
whence
\begin{align*}  
\det \big(Q_t ^{-1} -  Q_\infty^{-1}\big) & 
 = (\det Q_t)^{-1}\, e^{2t\,\text{tr} B}.
 \end{align*}
 Inserting this in \eqref{conto_contrvv}, we finally have
\begin{align}\label{conto_contr}
 I_t (u)
 &=
 (2\pi)^{-\frac n2}\,\int_{\R^{n}} \exp \left[
{-\frac12
\left\langle x' \,, x'\right\rangle}\right] \,
dx'=1.
\end{align}
This proves \eqref{th_Itu} and also that  $\gamma_{-\infty}$ is an invariant  measure for $\mathcal{H}_{t}^{UO}$.
Since $\left|\mathcal{H}_{t}^{UO} f\right| \le \mathcal{H}_{t}^{UO} |f|$, it also follows that each $\mathcal{H}_{t}^{UO} $ is a contraction on $L^1(\gamma_{-\infty})$. From \eqref{cons_eq} we see that $\mathcal{H}_{t}^{UO} $ is  a contraction also on $L^\infty(\gamma_{-\infty})$, and then by interpolation on $L^p(\gamma_{-\infty})$ for all $1\le p \le \infty$.
\end{proof}

From the proof of Proposition  \ref{pr:contr_L1_pd} we can also deduce the following formula, which may be of independent interest.
\begin{corollary}\label{semigroup_eRu}
For all $t>0$ and $x \in \R^n $ one has
\begin{equation}
\mathcal H_t^{OU}
\big( e^{R(\cdot)}\big) (x)=
e^{-t\,\mathrm{tr}B}\,
 e^{R(x)}. 
\end{equation}
\end{corollary}
\begin{proof}
 From \eqref{th_Itu} and \eqref{f:intermediate} we get
 
\begin{align}
1 =I_t (u)
&= c_{n}\,  e^{t\,\mathrm{tr} \,B}\,  e^{-R(u)}
\int_{\R^{n}}    \,K_t^{OU}(u,x)  \, 
\,dx
\notag\\
&=
 e^{t\,\mathrm{tr} \,B}\,  e^{-R(u)}
\int_{\R^{n}}    \,K_t^{OU}(u,x)  \, 
e^{R(x)}\, d\gamma_\infty (x)
\notag\\
&= 
 e^{t\,\mathrm{tr} \,B}\,  e^{-R(u)}
   \,\mathcal H_t^{OU}(
e^{R})(u).
\notag
\end{align}
The corollary is proved.
\end{proof}

\begin{proposition}\label{prop:strongcontinuity}
Suppose $p\in [1,\infty)$. Then $(\mathcal H_t^{UO})_{t>0}$ is strongly continuous on $L^{p}(\gamma_{-\infty})$.
\end{proposition}

\begin{proof}
It is known from  Proposition \ref{pr:contr_L1_pd} that
$\mathcal H_t^{UO}$ is bounded on $L^p(\gamma_{-\infty})$.

Suppose now $f\in C_{c}^{\infty}(\mathbb{R}^n)$. To prove that $\| \mathcal{H}_t^{UO} f - f \|_{L^{p}(\gamma_{-\infty})} \to 0$ as $t\to 0^{+}$, 
we shall apply the Dominated Convergence Theorem. 
Since $f$ is continuous and bounded, $\mathcal{H}_t^{UO} f(x)$ tends to $ f(x)$ pointwise as $t \to 0^+$. 
Therefore, it suffices to verify that there exists $\delta \in (0,1)$ such that the dominating function
\begin{equation*}
   \sup_{0<t<\delta} \left| \mathcal H_t^{UO}f(x) - f(x)\right|
\end{equation*}
belongs to $L^p(\gamma_{-\infty})$ for $1\le p < \infty$.

If $p > 1$, we can bound this supremum  by
\begin{equation*}
   \sup_{t>0} \left| \mathcal H_t^{UO}f(x) \right| + |f(x)| = \mathcal{H}_*^{UO}f(x) + |f(x)|,
\end{equation*}
where $\mathcal{H}_*^{UO}$ denotes the inverse Ornstein--Uhlenbeck maximal operator. Since $\mathcal{H}_*^{UO}$ is bounded on $L^p(\gamma_{-\infty})$ for all $1 < p \le \infty$ (see \cite[Theorem 3.1]{BCCS}), and $f \in C_c^\infty(\mathbb{R}^n) \subset L^p(\gamma_{-\infty})$, the dominating function
 belongs to $L^p(\gamma_{-\infty})$.
 
Assume now that $p = 1$. Using \eqref{litet}, \eqref{cons_eq} and the fact that $|D_{-t}\,x-u| \simeq |x-D_t\,u|$, one has
\begin{equation*}
    \left| \mathcal H_t^{UO}f(x) - f(x)\right| \lesssim e^{-R(x)}\, t^{-n/2}\,\int_{\mathbb{R}^n} \exp\left(-c\,\frac{|x-D_t\,u|^2}{t}\right)\, | f(u)- f(x)|\, e^{R(u)}\, du.
\end{equation*}
for $0<t<1$. Here we can estimate, applying \eqref{close} for the last step,
\begin{equation*}
   |x-D_t\,u|^2 =  |x- u + u - D_t\,u|^2 \ge  |x- u|^2 - 2 |x- u| |u - D_t\,u| \ge |x- u|^2 - C\, |x- u||u|t.
\end{equation*}

First, let us take $|x|>1+2\rho$, where $\rho = \sup\{ |u|\colon u \in \mathrm{supp}\, f\}$.
Then $f(x)=0$, it suffices to integrate over $u \in \mathrm{supp}\, f$, and we have $1 < |x-u| \simeq |x|$. If we choose $\delta$ and thus $t$ small enough, we obtain
\begin{equation*}
     t^{-n/2}\, \exp\left(-c\,\frac{|x-D_t\,u|^2}{t}\right)  \le   t^{-n/2}\,   \exp\left(-c\,\frac{|x-u|^2}{t} + C\, |x- u| |u|\right) \lesssim  \exp(-c\,|x|^2),
\end{equation*}
and therefore
\begin{equation*}
   \sup_{0<t<\delta} \left| \mathcal H_t^{UO}f(x) - f(x)\right| \lesssim e^{-R(x)}\, \exp(-c\,|x|^2) \int_{\R^{n}} | f(u)|\, e^{R(u)}\, du,
\end{equation*}
which is clearly in $L^1(\gamma_{-\infty})$.

Next, let $|x|<1+2\rho$ and consider the integral restricted to $|u|<2+4\rho$. Then this part of $\left| \mathcal H_t^{UO}f(x) - f(x)\right|$ is at most a constant times
\begin{equation*}
     t^{-n/2}\,\int_{|u|<2+4\rho} \exp\left(-c\,\frac{|x-u|^2}{t} + C \right) \, | f(u)- f(x)|\,  du
   \lesssim  \int_{|u|<2+4\rho} |x-u|^{1-n}\,  du \le C,
\end{equation*}
and the supremum in $t$ is thus integrable over  $|x|<1+2\rho$.

Finally, consider the case where $|x|<1+2\rho$ and $|u|>2+4\rho$, so that $|x-u| \simeq |u|$. Then
\begin{equation*}
     t^{-n/2}\, \exp\left(-c\,\frac{|x-D_tu|^2}{t}\right)  \le   t^{-n/2}\,  \exp\left(-c\,\frac{|x-u|^2}{t} + C\, |x- u||u| \right).
\end{equation*}
With $\delta$ small enough, this will be at most $C\, \exp(-c\,{|u|^2}/{t})$. Moreover, since $f$ is bounded, we can achieve
\begin{equation*}
    \int_{|u|>2+4\rho} \exp\left(-c\,\frac{|u|^2}{t} \right) \, | f(x)|\,  e^{R(u)} du \le C \sup |f|.
\end{equation*}
The last quantity is thus a uniform bound for the corresponding part of $\left| \mathcal H_t^{UO}f(x) - f(x)\right|$,  
which  ends the argument
 for $p=1$.

By the Dominated Convergence Theorem, $\| \mathcal H_t^{UO} f -f \|_{L^{p}(\gamma_{-\infty})}  \to 0$
 for all $1 \le p < \infty$ and $f \in C_c^\infty(\mathbb{R}^n)$. 
 The proposition now follows from this dense convergence, together 
 with the $L^{p}(\gamma_{-\infty})$ contractivity of $\mathcal H_{t}^{UO}$,
 via a standard density argument.                
  \end{proof}

\section{The spectrum of $\mathcal{A}^{Q,B}$}\label{section:spectrum}
By Proposition~\ref{prop:strongcontinuity}, for $1\le p<\infty$
the semigroup $(\mathcal H_t^{UO})_{t>0}$ has an infinitesimal
generator $\mathcal A_p^{Q,B}$ on $L^p(\gamma_{-\infty})$. We write $\sigma_p(\mathcal A^{Q,B})$ for its spectrum on
$L^p(\gamma_{-\infty})$, and generally omit the subscript $p$ when the underlying space is clear.

We start by proving a lemma.

\begin{lemma}\label{signM}
The eigenvalues of the matrix $M= -(B+QQ_\infty^{-1})$ have  negative real parts, $\mathrm{tr}\,M = \mathrm{tr}\,B$, and
\begin{equation}\label{samemeasure}
\int_0^\infty e^{sM}Qe^{sM^*}\, ds = Q_\infty.
\end{equation}
\end{lemma}

Notice that this means that $Q$ and $M$ satisfy the assumptions (H1) and (H2).
Thus we can associate an Ornstein–Uhlenbeck theory with the
operator  $\mathcal L^{Q,M}$,  and      
 \eqref{samemeasure} 
 implies that the corresponding invariant measure $\gamma_\infty$ is the same as that for $\mathcal L^{Q,B}$.         
 
\begin{proof}
Let us recall from~\cite[Lemma 2.1]{MPRS} that 
\begin{align}\label{BQinfBstar}
BQ_\infty+Q=-Q_\infty B^*,
\end{align}
which yields
\begin{align}\label{MBstarconj}
M = - (B+QQ_\infty^{-1})= Q_\infty B^*Q_\infty^{-1}.
\end{align}
Equation \eqref{MBstarconj} says that $M$ is similar to $B^*$. 
Therefore, the eigenvalues of $M$ are the complex conjugates of the eigenvalues of $B$, and thus the same as those of  $B$, with negative real parts.
Moreover, since $\mathrm{tr} (XY) = \mathrm{tr} (YX)$ 
and $B$ is a real matrix,
 it immediately follows that $\mathrm{tr} \, M = \mathrm{tr} \, B^* = \mathrm{tr} \, B$.

 By~\eqref{MBstarconj} moreover
\[
M^* = Q_\infty^{-1} B Q_\infty,
\]
and 
\[
e^{sM} = Q_\infty e^{sB^*} Q_\infty^{-1} \qquad \text{and} \qquad e^{sM^*} = Q_\infty^{-1} e^{sB} Q_\infty.
\]
Therefore,
\begin{align*}
\int_0^\infty e^{sM}Qe^{sM^*}ds
&= \int_0^\infty \left( Q_\infty e^{sB^*} Q_\infty^{-1} \right) Q \left( Q_\infty^{-1} e^{sB} Q_\infty \right) \, ds \\
  &= Q_\infty \left( \int_0^\infty e^{sB^*} (Q_\infty^{-1} Q Q_\infty^{-1}) e^{sB} \, ds \right) Q_\infty.
\end{align*}
By~\eqref{BQinfBstar} we have
\[
Q_\infty^{-1} Q Q_\infty^{-1} = Q_\infty^{-1} (-BQ_\infty - Q_\infty B^*) Q_\infty^{-1} = -Q_\infty^{-1}B - B^*Q_\infty^{-1}.
\]
Thus,
\begin{equation*}
\int_0^\infty e^{sM}Qe^{sM^*}ds
= Q_\infty \left( \int_0^\infty e^{sB^*} (-Q_\infty^{-1}B - B^*Q_\infty^{-1}) e^{sB} \, ds \right) Q_\infty.
\end{equation*}
It remains to observe that the integrand here is  the derivative of the function $s\mapsto -e^{sB^*} Q_\infty^{-1} e^{sB}$,
so that
\begin{equation*}
\int_0^\infty e^{sM}Qe^{sM^*}\, ds = Q_\infty \left[ -e^{sB^*} Q_\infty^{-1} e^{sB} \right]_0^\infty Q_\infty =  Q_{\infty},
\end{equation*}
and the proof is complete.
\end{proof}

\begin{lemma}\label{lem:adjoint_relation}
Let $U\colon L^2(\gamma_{\infty})  \to L^2(\gamma_{-\infty})$ be the isometric isomorphism given by multiplication by the factor 
$c_n^{1/2} e^{-R} $. For $ f\in L^{2}(\gamma_{-\infty})$, one has
\[
\mathcal{H}_t^{UO} = e^{t\,\mathrm{tr} \,B}\,    U\, ( \mathcal H_t^{OU})^* U ^{-1}, \qquad t>0.
\]
\end{lemma}

Here $(  \mathcal H_t^{OU})^*$ denotes the adjoint of $\mathcal{H}_t^{OU} $ in $L^2(\gamma_{\infty})$.

\begin{proof}
By \eqref{duenuclei} we have that
\begin{equation}\label{duenucleivv}
K_t^{UO}(x,u) =
c_{n}\,e^{-R(u)} \, e^{-R(x)}\, e^{t\,\mathrm{tr} \,B}\,K_t^{OU}(u,x).
\end{equation}
Thus
 \begin{align*}
 \mathcal H_t^{UO} f(x)&=\int_{\R^{n}} K_t^{UO}(x,u)\,f(u)\, d\gamma_{-\infty}(u)\\
&= c_{n}\, e^{t\,\mathrm{tr} \,B}\,
  e^{-R(x)}\,\int_{\R^{n}}  e^{-R(u)}\,K_t^{OU}(u,x)\,f(u)\, d\gamma_{-\infty}(u)\\
&= e^{t\,\mathrm{tr} \,B}\,
  e^{-R(x)}\,\int_{\R^{n}}  \,K_t^{OU}(u,x)\,f(u)\,  e^{R(u)}\, d\gamma_{\infty}(u),
\end{align*}
where we used \eqref{eq:rel_meas}.
Since  the kernel
 of the adjoint of the semigroup operator $ \mathcal H_t^{OU}$  with respect to the measure $\gamma_{\infty}$ is
 $ (K_t^{OU})^* (x,u) = K_t^{OU} (u,x)$, we conclude that
  \begin{align*}
 \mathcal H_t^{UO} f(x)
&= e^{t\,\mathrm{tr} \,B}\,
  e^{-R(x)}\,\int_{\R^{n}}  \, (K_t^{OU})^* (x,u)\,f(u)\,  e^{R(u)}\, d\gamma_{\infty}(u)\\
&=
 e^{t\,\mathrm{tr} \,B}\,    U\,
 \big(  \mathcal H_t^{OU}\big)^* U^{-1}f(x),
  \end{align*}
  which completes the proof.
\end{proof}

\begin{remark}\label{tracce}
 It follows from \eqref{MBstarconj} that
\[
 \mathrm{tr}  (QQ_\infty^{-1})= -(\mathrm{tr} \, M+\mathrm{tr} \, B)=-2\mathrm{tr} \, B.
\]
\end{remark}

The main result of this section is the following.

\begin{theorem}\label{spettro di A}
\begin{enumerate}
\item The spectrum $\sigma_{2}(\mathcal  A^{Q,B})$ of $\mathcal  A^{Q,B}$ in $L^2(\gamma_{-\infty})$ consists only of eigenvalues given by
\begin{align*}
\mu=\mathrm{tr} \, B-\sum_{j=1}^nm_j\rho_j,
\end{align*}
where $\rho_1,\ldots,\rho_n$
 are  the eigenvalues of the matrix $B+QQ^{-1}_\infty$ (counted with algebraic multiplicity)
 and the numbers $m_j$ are any non-negative integers.
In particular, $\mathcal  A^{Q,B}$ has a spectral gap.  Each of the corresponding eigenfunctions is of the form $e^{-R(x)}P(x)$, where $P$ is a polynomial.
\item $\mu \in \sigma_{2}(\mathcal{A}^{Q,B})$ if and only if there exists a {nonzero} homogeneous polynomial $P$ such that 
\[
\langle (B +QQ_{\infty}^{-1})x,\nabla P\rangle = (\mathrm{tr} \, B -\mu)P. 
\]
\end{enumerate}
\end{theorem}

 The proof of this theorem will be reduced to
 the study of an appropriate Ornstein--Uhlenbeck operator, 
whose spectrum is known. We shall use the following statement, which summarizes Proposition 3.1 and Lemma 3.2 in~\cite{MPP}.
\begin{proposition}\label{LemmaOUM}
Let $M$ be the matrix of Lemma~\ref{signM}, and consider
$\mathcal L^{Q,M}$ on $L^2(\gamma_\infty)$.
Let $\mu_1,\ldots,\mu_n$ be the eigenvalues of $M$, counted with multiplicity. Then
\begin{enumerate}
\item $\mathcal{L}^{Q, M}$ has a pure point spectrum.
Its eigenfunctions are polynomials;
the eigenvalues are given by the numbers $\sum_{j=1}^nm_j\mu_j$,
where the $m_j$
are nonnegative integers.
\item
A complex number $\lambda$ is an eigenvalue of $\mathcal{L}^{Q,M}$ if tnenand only if there is a {nonzero} homogeneous polynomial $P$ satisfying
\begin{align}\label{Mxcdotnabla}
\langle Mx, \nabla P\rangle =\lambda P.
\end{align}
\end{enumerate}
\end{proposition}
Before  the proof of Theorem~\ref{spettro di A}, we make some computations. We start by  noticing that  \cite[(4.1)]{CCS2} implies
\begin{align}\label{(4.1)}
2\langle Q^{-1}_\infty x, Bx\rangle =2 \langle B^*Q^{-1}_\infty x, x\rangle =-|Q^{1/2}Q^{-1}_\infty x|^2,
\end{align}
Moreover, 
we have
\begin{align}\label{grad gamma}
\nabla\gamma_{\infty}(x)= {c_{n}}\nabla e^{-\frac12 \langle Q^{-1}_\infty x, x\rangle }
=- c_{n}Q^{-1}_\infty x e^{-\frac12 \langle Q^{-1}_\infty x, x\rangle }
=-\gamma_{\infty}(x)Q^{-1}_\infty x
\end{align}
and because of \eqref{samemeasure} and \eqref{(4.1)}
\begin{align*}
\mathrm{tr}
\big( Q\nabla^2 \gamma_{\infty}\big)(x)
&=
-\mathrm{tr} (Q Q^{-1}_\infty)\gamma_{\infty}(x)
+\langle Q Q^{-1}_\infty x, Q^{-1}_\infty x\rangle \gamma_{\infty}(x)
\\
&=
(-\mathrm{tr} (Q Q^{-1}_\infty)+|Q^{1/2}Q^{-1}_\infty x|^2)\gamma_{\infty}(x).
\end{align*}
Combining  the last two identities and applying again \eqref{(4.1)}, we can compute the action of $\mathcal{A}^{Q,B}$ on $\gamma_{\infty}$,
that is,
\begin{align*}
\mathcal{A}^{Q,B}\gamma_{\infty}(x)
&= \frac12 \mathrm{tr} \big( Q\nabla^2 \gamma_{\infty}\big)(x) - \langle Bx, \nabla\gamma_{\infty}(x)\rangle
\\
&=
\frac12\left(|Q^{1/2}Q^{-1}_\infty x|^2-\mathrm{tr} (QQ^{-1}_\infty)\right)\gamma_{\infty}(x)
+ \langle Bx,  Q^{-1}_\infty x\rangle \gamma_{\infty}(x)
\\
&= -\frac12 \mathrm{tr} (Q Q^{-1}_\infty)\gamma_{\infty}(x).
\end{align*}
Recalling that $\mathrm{tr} (QQ^{-1}_\infty) = -2\mathrm{tr} \, B$ (see Remark \ref{tracce}), we immediately conclude that
\begin{align}\label{autofunzione gamma}
\mathcal{A}^{Q,B}\gamma_{\infty}=\mathrm{tr} \, B\,\gamma_{\infty}.
\end{align}

\begin{proof}[Proof of Theorem \ref{spettro di A}]
The formulae  \eqref{autofunzione gamma}  and \eqref{grad gamma}  yield
\begin{equation}\label{APgamma}
\begin{split}
\mathcal{A}^{Q,B}(\gamma_{\infty} f)
&=\gamma_{\infty} \mathcal{A}^{Q,B}f+f\mathcal{A}^{Q,B}\gamma_{\infty}+\langle Q \nabla\gamma_{\infty}, \nabla f\rangle 
\\
&=\gamma_{\infty} \mathcal{A}^{Q,B} f+\mathrm{tr} \, B\gamma_{\infty} f-\gamma_{\infty} \langle QQ_\infty^{-1}x, \nabla f\rangle
\\
&=\frac12\gamma_{\infty} \mathrm{tr} (Q\nabla^{2}f)+\mathrm{tr} \, B \gamma_{\infty} f-\gamma_{\infty} \langle (B+Q Q_\infty^{-1})x, \nabla f\rangle,
\end{split}
\end{equation}
and hence
\begin{equation}\label{unitaryeq}
\gamma_{\infty}^{-1} \mathcal{A}^{Q,B}(\gamma_{\infty} f) = \frac12  \mathrm{tr} (Q\nabla^{2}f) - \langle (B+QQ_\infty^{-1})x, \nabla f\rangle +\mathrm{tr} \, B f.
\end{equation}
If $U\colon L^2(\gamma_{\infty}) \to L^2(\gamma_{-\infty})$ is the isometry of Lemma~\ref{lem:adjoint_relation}, setting
as before
\[
M = - (B+QQ_\infty^{-1}),
\]
then we may rewrite \eqref{unitaryeq}  as
\begin{equation}\label{unitaryeq2}
\mathcal{A}^{Q,B} f = U \big( \mathcal{L}^{Q,M} + \mathrm{tr} \, B \big) U^{-1} f.
\end{equation}
Since $C_c^\infty(\mathbb R^n)$ is a core for
$\mathcal L^{Q,M}$ on $L^2(\gamma_\infty)$ and $U$ maps $C_c^\infty(\mathbb R^n)$
onto itself, the above is an equality also of domains. Thus, the operator $\mathcal{A}^{Q,B}$ in  $L^{2}(\gamma_{-\infty})$ is unitarily equivalent to $\mathcal{L}^{Q,M} + \mathrm{tr} \, B$ in  $L^{2}(\gamma_{\infty})$.

This implies that
\[
\sigma_{2}(\mathcal{A}^{Q,B}) = \mathrm{tr} \, B + \sigma_{2}(\mathcal{L}^{Q,M}),
\]
where $\sigma_{2}(\mathcal{L}^{Q,M})$ 
denotes the spectrum of $\mathcal{L}^{Q,M}$ on $L^{2}(\gamma_{\infty})$. 
The theorem follows by Proposition~\ref{LemmaOUM}.
\end{proof}

\begin{remark}\label{remark:orthogonality} 
As a consequence of the unitary equivalence between 
$\mathcal{A}^{Q,B}$ and $\mathcal{L}^{Q,M} + \mathrm{tr} \, B$, 
we can deduce some orthogonality properties
  for the generalized    eigenspaces of $\mathcal{A}^{Q,B}$ 
  from the   Ornstein--Uhlenbeck setting. 
By \cite[Proposition 3.1]{CCS4}, 
if the  matrix $M$ has only   one eigenvalue, 
then the generalized eigenspaces of $\mathcal{L}^{Q,M}$ 
are mutually orthogonal in $L^2(\gamma_{\infty})$. 
Since $M$ is 
similar to $B^*$ (see Lemma \ref{signM}), 
$M$ has a unique eigenvalue if and only if $B$ has.  Therefore, if $B$ has only one  eigenvalue, any two distinct generalized eigenspaces 
of $\mathcal{A}^{Q,B}$ are orthogonal in $L^2(\gamma_{-\infty})$. 
But  if $B$ has distinct eigenvalues, orthogonality may fail, 
as  can be seen from an example given in  \cite[Section 5]{CCS4}.
\end{remark}

 \begin{proposition} 
\begin{enumerate}
 \item If $p\in (1,\infty)$, then $\sigma_{p}(\mathcal{A}^{Q,B}) = \sigma_{2}(\mathcal{A}^{Q,B}) $,
  and each point of   $\sigma_{p}(\mathcal{A}^{Q,B})$    is an eigenvalue.            
 \item If $p=1$, then $\sigma_{1}(\mathcal{A}^{Q,B}) = \{z\in \mathbb{C} \colon \mathrm{Re}\, z \leq 0 \}  $. 
 \end{enumerate}
\end{proposition}

\begin{proof}

 Let us first prove (1), thus assuming $1<p<\infty$.

We start with the Ornstein--Uhlenbeck semigroup $\mathcal{H}_t^{OU}$ 
 in $L^2(\gamma_{\infty})$. From \cite[Section 3(c)]{Lunardi} we have that $\mathcal{H}_t^{OU}$ is a compact operator on $L^2(\gamma_{\infty})$ for each $t>0$.
The resolvent of $\mathcal L^{Q,B}$ satisfies
$$
\left(\lambda - \mathcal L^{Q,B}  \right)^{-1} = \int_{0}^{\infty} e^{-t\lambda}\,\mathcal{H}_t^{OU}\,dt
$$
for $\lambda>0$.  Here, and similarly in the sequel, we write simply $\lambda$ for $\lambda I$. This resolvent                                           
is also compact for  $\lambda>0$ (see \cite[p.~41]{MPP}  or \cite[p.~12]{Lunardi}), and by Lemma~\ref{signM} the same property holds for $(\lambda - \mathcal{L}^{Q,M})^{-1}$, where $M = -(B+QQ_\infty^{-1})$.
Further,  the nonzero part of the spectrum of  $\left(\lambda - \mathcal L^{Q,B}  \right)^{-1}$   consists only of eigenvalues.

\medskip

We then apply \cite[Ch.\ II, Theorem 1.10, p.~55]{EngelNagel}, 
with $w =0$, $M=1$ and $\lambda$ replaced by $\lambda + \mathrm{tr}\,B$, 
to the contractive semigroup $\left(\mathcal{H}_t^{UO}\right)_{t>0}$ 
in the Banach space $L^p(\gamma_{-\infty})$, to get that 
for ${\Re}\,\lambda + \mathrm{tr}\,B > 0$ 
the resolvent $(\lambda - \mathcal A^{Q,B} + \mathrm{tr}\,B)^{-1}$ 
is a bounded operator on $L^p(\gamma_{-\infty})$, given by
$$(\lambda - \mathcal A^{Q,B} + \mathrm{tr}\,B )^{-1}\,f 
= \int_0^\infty e^{-(\lambda + \mathrm{tr}\,B )t}\,\mathcal{H}_t^{UO}\,f \, dt, \qquad f \in L^p(\gamma_{-\infty}).$$

We now fix a real $\lambda > -\mathrm{tr}\,B > 0$. 
With the isometric isomorphism $U: L^2(\gamma_{\infty}) \to L^2(\gamma_{-\infty})$ 
 defined above, \eqref{unitaryeq2} says that $\mathcal A^{Q,B} - \mathrm {tr}B= U \mathcal L^{Q,M} U^{-1}$, and it follows that
$$
\left(\lambda - \mathcal A^{Q,B} + \mathrm {tr}B \right)^{-1} = U (\lambda - \mathcal L^{Q,M})^{-1} U^{-1}.
$$
Since $(\lambda - \mathcal{L}^{Q,M})^{-1}$ is compact on $L^2(\gamma_\infty)$, 
it follows that $(\lambda - \mathcal{A}^{Q,B} + \mathrm{tr}\,B)^{-1}$ is compact on $L^2(\gamma_{-\infty})$,  and its nonzero spectrum contains
   only eigenvalues.

By~\cite[Theorem 1.6.1 and Corollary 1.6.2]{Davies},  
 the compactness of $(\lambda - \mathcal{A}^{Q,B} + \mathrm{tr}\,B)^{-1}$ on $L^2(\gamma_{-\infty})$, 
 combined with its boundedness on $L^p(\gamma_{-\infty})$ for all $1 < p < \infty$, 
 ensures that $(\lambda - \mathcal{A}^{Q,B} + \mathrm{tr}\,B)^{-1}$ 
 is compact on $L^p(\gamma_{-\infty})$ and that its spectrum $\sigma_p\big((\lambda - \mathcal{A}^{Q,B} + \mathrm{tr}\,B)^{-1}\big)$ 
 is independent of $p \in (1,\infty)$.  Notice   that Davies' arguments work also under the modified hypothesis $1<p_1<p<p_0<\infty$.

To pass from the spectrum of $(\lambda  - \mathcal{A}^{Q,B} + \mathrm {tr}B)^{-1}$ to that of  $\mathcal{A}^{Q,B}$, we  first observe that, by
 \cite[Ch.\ II, Theorem 1.4, p. 51]{EngelNagel},
 the generator $\mathcal{A}^{Q,B}$ of the semigroup  $\left(\mathcal{H}_t^{UO}\right)_{t>0}$ 
  is a closed operator in $ L^p(\gamma_{-\infty})$. This allows us to apply  the Spectral Mapping Theorem, see \cite[Ch.\ IV, Theorem 1.13(i), p.\ 243]{EngelNagel}. As a result,
\begin{equation} \label{fff}
 \sigma_p\left((\lambda  - \mathcal{A}^{Q,B} + \mathrm {tr}B)^{-1}\right) \setminus\{0\} 
 = \left\{  \frac 1{\lambda + \mathrm {tr}B - \mu}:\, \mu \in \sigma_p(\mathcal{A}^{Q,B})\right\},
  \end{equation}
 or equivalently,
 \begin{equation} \label{ett}
  \sigma_p(\mathcal{A}^{Q,B})  
   = \left\{  \lambda + \mathrm {tr}B -\frac 1 \nu: \, \nu \in \sigma_p\left((\lambda - \mathcal{A}^{Q,B} + \mathrm {tr}B)^{-1}\right) \setminus\{0\}\right\}. 
  \end{equation} 
 Since the right-hand side in \eqref{ett} 
 is independent of $p \in (1, \infty)$, so is $\sigma_p(\mathcal{A}^{Q,B})$.           

 Part (ii) of the same theorem \cite[Ch.~IV, Theorem~1.13, p.~243]{EngelNagel} says that \eqref{fff} and thus \eqref{ett}
hold also for the point spectra.
Since $(\lambda - \mathcal{A}^{Q,B} + \mathrm{tr}\,B)^{-1}$ is compact on $L^p(\gamma_{-\infty})$, the set $\sigma_p\big((\lambda - \mathcal{A}^{Q,B} + \mathrm{tr}\,B)^{-1}\big)\setminus\{0\}$ coincides with the point spectrum of this resolvent, and hence the same holds for 
$\mathcal{A}^{Q,B}$. This ends the proof of (1).

Let us now prove (2). As in \eqref{MBstarconj},
we write $M = - (B+QQ_\infty^{-1})$.
The map
\[
\widetilde U = c_n^{-1/2} U \colon L^1(\mathbb{R}^n, dx) \to L^1(\gamma_{-\infty}), \qquad \widetilde Uf :=e^{-R(x)}  f,
\]
is an isometric isomorphism. 
Since
$C_c^\infty(\mathbb R^n)$ is a core for the
Ornstein--Uhlenbeck generator on $L^1(dx)$
by \cite[Proposition 3.2]{Metafune},
and $\widetilde U$ maps $C_c^\infty(\mathbb R^n)$ onto itself, the spectrum of $\mathcal{A}^{Q,B}$ on $L^1(\gamma_{-\infty})$ 
coincides with that of $ \mathcal{L}^{Q,M} + \mathrm{tr} \, B$ on $L^1(\mathbb{R}^n, dx)$. 
By \cite[Theorem 4.12]{Metafune}, the spectrum of $\mathcal{L}^{Q,M}$ on $L^1(\mathbb{R}^n, dx)$
  is the half-plane $\{z \in \mathbb{C} \colon \mathrm{Re}\, z \leq -\mathrm{tr} \, M\}$. 
 Recalling from Lemma \ref{signM} that 
 $\mathrm{tr} \, M = \mathrm{tr} \, B$, 
 we conclude that the spectrum of $\mathcal{A}^{Q,B}$ on $L^1(\gamma_{-\infty})$ 
 is the half-plane $\{z \in \mathbb{C} \colon \mathrm{Re}\, z \leq 0\}$.
\end{proof}

 \begin{remark}
The domain $D(\mathcal{A}^{Q,B})$ of the infinitesimal generator in $L^2(\gamma_{-\infty})$ can be characterized by means of the isometry $U$ introduced above. More precisely, for $p=2$, the unitary equivalence  \eqref{unitaryeq2}          
implies that
\[
D(\mathcal{A}^{Q,B})_{L^2(\gamma_{-\infty})} = \{ f \in L^2(\gamma_{-\infty}) : U^{-1}f \in W^{2,2}(\gamma_{\infty}) \},
\]
where $W^{2,2}(\gamma_{\infty})$ is the domain of the classic Ornstein--Uhlenbeck operator $\mathcal{L}^{Q,M}$ in $L^2(\gamma_{\infty})$ \cite[Section 3]{MPRS}.
\end{remark}

\section{Ultracontractivity of  $(\mathcal H_t^{UO})_{t>0}$}\label{section:ultra}

We recall that a semigroup  $(T_t)_{t>0}$
is called ultracontractive if for every $p\ge 1$ and $t >0$ there exists some positive constant
$c_{p}(t)$
 such that
 $\|T_t f\|_{\infty}\le  c_{p}(t)
\|f\|_{p}$.

\medskip

It follows from earlier results by Nelson \cite{Nelson}  that the Ornstein--Uhlenbeck semigroup $(\mathcal H_t^{OU})_{t>0}$
is not 
ultracontractive in $L^p(\gamma_\infty)$.
 But $(\mathcal H_t^{OU})_{t>0}$
is ultracontractive in $L^p$
spaces endowed with Lebesgue measure { \cite[p.~13]{Lunardi}}. 
Unexpectedly, in the inverse Gaussian setting it is easy to prove ultracontractivity for $(\mathcal H_t^{UO})_{t>0}$
in $L^p$ spaces both with respect to Lebesgue measure
and with respect to $\gamma_{-\infty}$.

\smallskip

We shall need an elementary lemma.
\begin{lemma}\label{l:lr}
For $1\le r \le \infty$
we have
\begin{align*}
 \| e^{-\frac12 \langle Q_t^{-1}(e^{tB}u-x) , \, e^{tB}u-x\rangle}\|_{L^{r}(du)}&
   = \left(\frac{2\pi} r\right)^{\frac{n}{2r}} \left(\det Q_t\right)^{\frac1{2r}}\,e^{-\frac{t}{r}\, \mathrm{tr}B}, \qquad x \in \R^n,
 \end{align*}
 where the right-hand side is $1$ for $r=\infty$.
\end{lemma}

\begin{proof}
When $r=\infty$ this is trivial, and for  $1\le r < \infty$ one has
\begin{align*}
\int_{\R^n} e^{-\frac{r}2 \langle Q_t^{-1/2}(e^{tB}u-x) ,\,  Q_t^{-1/2}(e^{tB}u-x)\rangle}du
= (\det Q_t)^{1/2}\, e^{-t\, \text{tr}B}\, \int_{\R^n} e^{-\frac{r}2  |z|^2 }\,dz,
\end{align*}
which implies the lemma.
 \end{proof}

\begin{proposition}
\begin{enumerate}
\item
The semigroup 
 $\big(\mathcal H_t^{UO}\big)_{t>0}$
 is ultracontractive with respect to Lebesgue measure
 with constant 
 \begin{equation}
    \label{def:cpt}
    \begin{split}
c_{p}(t) & = (2\pi)^{-\frac{n}{2p}}\, e^{t\,\mathrm{tr}B \frac1{p}} 
 {(p')}^{-\frac{n}{2p'}} \left(\det Q_t\right)^{-\frac1{2p}} ,  \qquad 1<p<\infty;  \\
 c_{1}(t) & = (2\pi)^{-\frac{n}{2}}\, e^{t\,\mathrm{tr}B } \left(\det Q_t\right)^{-\frac1{2}}.
    \end{split}
\end{equation}
\item
The semigroup 
 $\big(\mathcal H_t^{UO}\big)_{t>0}$
 is ultracontractive  with respect to $\gamma_{-\infty}$
with  the same constant 
 $c_{p}(t)$ given by \eqref{def:cpt}.
  \end{enumerate}
\end{proposition}
\begin{proof}
We may assume $f\in \mathcal C_c^{\infty}(\R^n)$ since this space is dense in 
$ L^p(du)$, for $1\leq p<\infty$.
Then it follows from \eqref{kolmogo} that
\begin{equation}\label{eq:kernel_integral}
\mathcal H_t^{UO}
f(x)=(2\pi)^{-\frac{n}{2}}
({\det} \, Q_t)^{-\frac{1}{2} }\, e^{t\,\text{tr}B}
\int_{\R^{n}} f(u)  \, e^{-\frac12 \langle Q_t^{-1}(e^{tB}u-x) , \, e^{tB}u-x\rangle}\, du\,, \qquad x\in\R^n,
\end{equation}
 and Lemma \ref{l:lr}  implies the assertion (1).

 To prove the assertion (2), it is enough to observe that
 \[\|f\|_{L^p(du)}\le \|f\|_{L^p(\gamma_{-\infty})},\]
 which completes the proof.
 \end{proof}

\begin{corollary}\label{cor:5.3}
 For any $ 1\le p<q$ and $t > 0$
\begin{equation}\label{transition}
\|\mathcal{H}_t^{UO} f\|_{L^q(\gamma_{-\infty})} \le \big(c_p(t)\big)^{1-\frac{p}{q}} \|f\|_{L^p(\gamma_{-\infty})}. 
\end{equation}
\end{corollary}
\begin{proof}
This  is immediate by interpolation  between the $L^p$-contractivity in Theorem \ref{pr:contr_L1_pd} and the ultracontractivity bound.
\end{proof}

\begin{remark}

Thanks to the explicit expression of the ultracontractivity constant $c_p(t)$ in \eqref{def:cpt}, we notice what follows.
For small $t $ one has 
\[
Q_t = \int_0^t e^{sB} Q e^{sB^*} ds \simeq t Q
\]
and consequently 
\[
\det Q_t \simeq t^n \det Q,
\]
which implies that the constant $c_p(t)$  blows up like $t^{-\frac{n}{2p}}$ as $t \to 0^+$.
Notice that the blow-up rate $t^{-\frac{n}{2p}}$ 
is analogous to that of the standard heat semigroup $(e^{t\Delta})_{t>0}$ in $\mathbb{R}^n$ 
from $L^p(dx)$ to $L^\infty(dx)$. This is a local effect. 
As $t \to +\infty$, the matrix $Q_t$ converges to the positive definite matrix $Q_\infty$. Since  $\textrm{tr}  B < 0$,
 the factor  $e^{\frac{t}{p}\textrm{tr}B}$ makes  $c_p(t)$  decay exponentially as
 $t \to +\infty$.

\end{remark}
\subsection{A counterexample for $q<p$}\label{count1}

 Clearly, \eqref{transition} provides bounded operators $L^p \to L^q$ for any $t$ as soon as $q  > p$.
   Let now instead $1<q<p<\infty$.


We claim that if  $t $ is small enough, in a way that  depends only on $p, q, n,\,Q$ and $B$, then $\mathcal H_t^{UO}$ is not bounded from $L^p(\gamma_{-\infty})$ to $L^q(\gamma_{-\infty})$.

Fix  $t \in (0,1)$ and let $u_0$ be a point in $\R^n$ with  $|u_0|> t^{-1/2}$. We   let $f$ be the characteristic function of the ball  $B(u_0, 1/|u_0|)$.   Since the density  $e^{R(u)}$ has constant order of magnitude in this ball, we have
\[
\|f \|_{L^p(\gamma_{-\infty})}  \simeq e^{R(u_0)/p}\,|u_0|^{-n/p}.
\]          
We will evaluate $\mathcal H_t^{UO} f(x)$ at points $x$ in the ball $B(D_t\,u_0,\sqrt t)$, and actually in a thin, curved slice of this ball, given by 
\[
S = \{x \in B(D_t\,u_0,\sqrt t): \, |R(x) - R(D_t\,u_0)| < 1\}, 
\] 
thus containing the center of the ball. 
Since $t<1$ we have $R(D_t\,u_0) \simeq R(u_0) \simeq |u_0|^2$. The width of $S$ is roughly $1/|u_0|$ as seen via $R(x) $, and its Lebesgue measure will be
\[
|S| \simeq t^{(n-1)/2}/|u_0|. 
\]
For $x \in S$ and  $u \in B(u_0, 1/|u_0|)$ we have 
$$
|D_{-t}\,x - u| \simeq |x - D_{t}\,u| \le |x - D_{t}\,u_0| +   |D_{t}\,(u - u_0)|    \lesssim       \sqrt t +  1/|u_0| \lesssim \sqrt t,
$$
the last estimate by the choice of $u_0$,
and \eqref{litet} implies that $ K_t^{UO}(x,u) \gtrsim t^{-n/2}\,e^{-R(x)} \simeq t^{-n/2}\,e^{-R(D_t\,u_0)} $. Thus
{\[
\mathcal H_t^{UO} f(x) \gtrsim  t^{-n/2}\,e^{-R(D_t\,u_0)} \, e^{R(u_0)}\, |u_{0}|^{-n}, \qquad x \in S,
\]  
 and so
\begin{align*}
\|\mathcal{H}_t^{UO}f\|_{L^q(\gamma_{-\infty})}
&\gtrsim t^{-n/2}e^{-R(D_tu_0)}e^{R(u_0)} |u_0|^{-n}\gamma_{-\infty}(S)^{1/q} \\
&\simeq t^{-n/2}e^{-R(D_tu_0)}e^{R(u_0)} |u_0|^{-n}|S|^{1/q}e^{R(D_tu_0)/q} \\
&\gtrsim e^{R(D_tu_0)(-1+1/q)}e^{R(u_0)}t^C|u_0|^{-C},
\end{align*}}
for some constants $C$.
For the quotient between the norms, we find
 \begin{align*}
     P := \|\mathcal H_t^{UO} f \|_{L^q(\gamma_{-\infty})}/\|f \|_{L^p(\gamma_{-\infty})} \gtrsim
e^{R(D_t\,u_0)(-1+1/q)}\,e^{R(u_0)(1-1/p)}\,t^{C}\,|u_0|^{-C}.
\end{align*}

In the light of \eqref {vel-4},  $\partial_s R(D_{s} \, u_0) \simeq |D_{s} \, u_0|^2 \simeq |u_0|^2$ for  $s \in (0,1)$, so that
\[
R(D_t\,u_0) - R(u_0) \simeq \int_{0}^{t} |u_0|^2\,ds = t|u_0|^2. 
\] 
Thus
\begin{align}\label{ineq:estP}
     P  \gtrsim e^{R(u_0)(1/q-1/p)}\, e^{-C t|u_0|^2}
\,t^{C}\,|u_0|^{-C}.  
\end{align}

 Recall that $R(u_0) =\frac12 \langle Q_\infty^{-1}u_0,u_0\rangle \simeq |u_0|^2$. By taking $t$ small in the way described in the claim, we can thus obtain the inequality $Ct|u_0|^2 < R(u_0)(1/q-1/p)/2$ between the first two exponents in \eqref{ineq:estP}. This leads to
\begin{align*}
     P  \gtrsim e^{R(u_0)(1/q-1/p)/2}\,t^{C}\, |u_0|^{-C},   
\end{align*}
which tends to infinity with $|u_0|$. The claim follows.

\section{Hypercontractivity}\label{sec:hyper}

We recall that a semigroup $(T_t)_{t>0}$
  is said to be hypercontractive if for every $p>1$ and 
$t>0$       there exists some $q(t)>p$
 such that
\begin{equation}\label{d:hyper}
 \|T_t f\|_{q(t)}\le 
\|f\|_{p}. 
\end{equation}

\begin{theorem}\label{hyper}
The semigroup $(\mathcal{H}_{t}^{UO})_{t>0}$ is hypercontractive.
\end{theorem}

\begin{proof}

 The main hypercontractivity property of $(\mathcal H_t^{OU})_{t>0}$ is the following.
Given $1\le p<q<\infty$ and $t>0$, if 
    \[
    q - 1 \leq  (p - 1) \|Q_\infty^{-1/2} e^{tB} Q_\infty^{1/2} \|^{- 2}
    \]
    then $\mathcal H_t^{OU}$ is a contraction from $L^p(\gamma_\infty)$ into $L^q (\gamma_\infty)$.
For a proof of this result,  holding also in the infinite-dimensional setting, we refer to  \cite{CM} and  \cite{MF}.

It has been proved in 
 \cite[Proposition 4.2 and Proposition 4.1]{MF} that 
\begin{equation}\label{in:MF}
 \|Q_\infty^{-1/2} e^{tB} Q_\infty^{1/2} \|  <1\qquad \text{ for all $t>0$.}
\end{equation}

We fix $t>0$.  Choosing $q = 2$ above, we see that we can find a $p  \in (1, 2)$ such that  $\mathcal{H}_t^{OU} $ is a contraction from $L^p(\gamma_{\infty})$ to $L^2(\gamma_{\infty})$. For the dual operator we then conclude that 
$\big(  \mathcal H_t^{OU}\big)^*$ is a contraction from  $L^2(\gamma_{\infty})$ to $L^{p’}(\gamma_{\infty})$, and thus from  $L^2(\gamma_{\infty})$ to any $L^{\varrho}(\gamma_{\infty})$ with $\varrho \in (2, p')$.

It follows that $\big(  \mathcal H_t^{OU}\big)^*\,U^{-1}$ is a contraction from  $L^2(\gamma_{-\infty})$ to $L^{\varrho}(\gamma_{\infty})$, for  $\varrho \in (2, p')$. Then we observe that
$U$ is a bounded operator   from  $L^{\varrho}(\gamma_{\infty})$ to $L^{\varrho}(\gamma_{-\infty})$, actually for any $\varrho \ge 2$; more precisely we have with $0 \le g \in  L^\varrho(\gamma_{\infty})$
 \begin{align*}
    \int_{\R^{n}} (Ug)^\varrho\,d\gamma_{-\infty} 
    &= c_n^{\varrho/2} \,\int_{\R^{n}} e^{-\varrho R(x)}\,g(x)^\varrho\,e^{R(x)}\,dx \\
    &\le  c_n^{\varrho/2} \,\int_{\R^{n}} g(x)^\varrho\,e^{-R(x)}\,dx =  c_n^{\varrho/2-1} \,\int_{\R^{n}} g^\varrho\,d\gamma_{\infty}.
\end{align*}
    Thus the operator norm is no larger than $c_n^{1/2-1/\varrho}$.

Now 
Lemma  \ref{lem:adjoint_relation} implies that $\mathcal{H}_t^{UO} $   is bounded from $ L^2(\gamma_{-\infty})$ to $  L^{\varrho}(\gamma_{-\infty})$
with operator norm at most $e^{t\,\mathrm{tr} \,B}\,c_n^{1/2-1/\varrho}$, if  $2 < \varrho < p'$.   To make it  a contraction, it is enough to take 
$\varrho \in (2, p')$ such that $1/2-1/\varrho \le t|\mathrm{tr} \,B|/\log_+ c_n $ (if $c_n \le 1$ any $\varrho \in (2, p')$ will do).

Interpolation with the contractive estimate on $L^1(\gamma_{-\infty})$ yields a contractive estimate $L^p(\gamma_{-\infty}) \to L^q(\gamma_{-\infty})$ for some $q > p$ for every $1 < p < 2$. Interpolation with the contractive estimate on $L^\infty$ similarly yields the stated contractivity for $p > 2$.
\end{proof}

As a consequence of the above lemma, we also have the following.
\begin{corollary}\label{rem:stability}
The
semigroup $(\mathcal{H}_t^{UO})_{t \ge 0}$ is exponentially stable on $L^2(\gamma_{-\infty})$, 
in the sense that it satisfies a bound           
\begin{equation}\label{eq:exp_bound}
    \|\mathcal{H}_t^{UO}\|_{L^2(\gamma_{-\infty}) \to L^2(\gamma_{-\infty})} \le e^{- \omega t} \qquad \text{for all } t \ge 0,
\end{equation}
with  $\omega = - \mathrm{tr} \,B > 0$.
\end{corollary}

\end{document}